\documentclass[reqno, 12pt]{amsart}
\usepackage{amssymb,amsmath,amsthm,newlfont,enumerate}

\renewcommand{\tilde}{\widetilde}

\theoremstyle{plain}
\newtheorem{theorem}{Theorem}[section]
\newtheorem{lemma}[theorem]{Lemma}

\theoremstyle{definition}

\theoremstyle{remark}

\newcommand{\CC}{{\mathbb C}}
\newcommand{\DD}{{\mathbb D}}

\newcommand{\ZZ}{{\mathbb Z}}

\newcommand{\cF}{{\mathcal F}}

\DeclareMathOperator{\hol}{Hol}

\DeclareMathOperator{\dist}{dist}
\DeclareMathOperator{\diam}{diam}

\begin{document}

\title[Landau constants]{Landau constants for algebroid multifunctions}

\author[T. Ransford]{Thomas Ransford}

\address{%
D\'epartement de math\'ematiques et de statistique, Universit\'e Laval, Qu\'ebec (QC), G1V 0A6, Canada}

\email{ransford@mat.ulaval.ca}

\thanks{Research supported by  NSERC Discovery Grant RGPIN-2026-04565}



\date{30 October 2026}

\begin{abstract}
In an earlier paper, Bloch's theorem was extended to multifunctions on the unit disk of the form
$$
F(z)=\{w\in\CC: w^n+a_1(z)w^{n-1}+\dots+a_{n-1}(z)w+a_n(z)=0\},
$$
where $n$ is an integer with $n\ge1$ and $a_1,\dots,a_n$ are holomorphic functions on the unit disk and $a_1'(0)=1$.
In this paper, we quantify this result by showing that the image of $F$ always contains a disk of radius 
$1/(2^{15n}n^{6n})$. Denoting by $L_n$ the supremum of all radii with this property (the Landau constant
of degree $n$), we further show that $L_n\le (3\pi\sqrt{3}/2)^{1/2}n^{-3/2}+O(n^{-2})$ as $n\to\infty$.
\end{abstract}

\maketitle

\section{Introduction}

Let $\DD$ denote the open unit disk in the complex plane $\CC$.
Given a domain $D$ in~$\CC$, we write $\hol(D)$ for the set of holomorphic functions on $D$.
Also, if $S$ is any subset of $\CC$,
then we write $\rho(S)$ for  its \emph{inradius}, namely
\[
\rho(S):=\sup\{r\ge0: \text{$S$ contains a disk of radius $r$}\}.
\]

The \emph{Landau constant} is defined as
\begin{equation}\label{E:Ldef}
L:=\inf\{\rho(f(\DD)): f\in\hol(\DD),\, f'(0)=1\}.
\end{equation}
It is a non-trivial fact that $L>0$. This is part of Bloch's theorem,
and is of interest, among other reasons, because it leads to a relatively simple proof of Picard's theorem for entire functions
(see e.g.\ \cite[Chapter~XII]{Co78}). The exact value of $L$ is not known, but it satisfies
\[
0.51 <L\le \frac{\Gamma(1/3)\Gamma(5/6)}{\Gamma(1/6)}\approx 0.54325\dots
\]
The upper bound is due to Rademacher \cite{Ra43}, and is
widely conjectured to be the true value of $L$. The lower
bound is from the recent preprint of Wikstr\"om \cite{Wi26}.

In \cite{GR05}, Guillot and the author established an analogue
of Bloch's theorem for algebroid multifunctions.
We recall here the terminology and the principal results.

An \emph{algebroid multifunction} of degree $n$ on a domain $D$ is a
set-valued function $F$ of the form
\begin{equation}\label{E:algdef}
F(z)=\Bigl\{w\in\CC: w^n+a_1(z)w^{n-1}+\cdots+a_{n-1}(z)w+a_n(z)=0\Bigr\},
\end{equation}
where $a_1,\dots,a_n\in\hol(D)$.
It is understood that the roots $w$ are counted according to multiplicity.
The \emph{image} of $F$ is defined as
\[
F(D):=\cup_{z\in D}F(z).
\]

Given a family $\cF$ of algebroid multifunctions of degree $n$
defined on the unit disk~$\DD$, we define the associated Landau constant by analogy
with \eqref{E:Ldef}:
\[
L(\cF):=\inf\{\rho(F(\DD)): F\in\cF\}.
\]
The goal of \cite{GR05}, and of the article \cite{LR00} that preceded it, was to establish  conditions on $\cF$ under which $L(\cF)>0$.
The main result proved was \cite[Theorem~1]{GR05},
to the effect that $L(\cF_n)>0$,
where $\cF_n$ is the set of $F$ as in \eqref{E:algdef} for which $a_1'(0)=1$.
This choice of normalization was discussed at length in \cite{GR05},
see in particular Theorem~5 of that paper.

As in \cite{GR05}, let us relabel $L(\cF_n)$ as $L_n$.
Then $L_1=L$, the classical Landau constant, and it seems
reasonable to call $L_n$ the Landau constant of degree~$n$. 
Our goal in this paper is to establish effective upper and lower bounds for $L_n$.

We begin with the lower bound.
The proof in \cite{GR05} that $L_n>0$
was by contradiction, so it yielded no quantitative
lower bound for $L_n$. Our first result furnishes such a bound.
 Moreover, in contrast with the techniques used in \cite{GR05} and \cite{LR00},
the proof is entirely elementary.

\begin{theorem}\label{T:lb}
For all $n\ge1$, we have 
\begin{equation}\label{E:lb}
L_n>\frac{1}{2^{15n}n^{6n}}.
\end{equation}
\end{theorem}

As for upper bounds,
we have $L_n\le L/n$ for all $n$.
Indeed, it was shown in \cite[Theorem~2]{GR05} that 
\begin{equation}\label{E:subadd}
\frac{1}{L_{n+m}}\ge \frac{1}{L_n}+\frac{1}{L_m} \quad(n,m\ge1),
\end{equation}
from which the bound $L/n$ readily follows.
In fact,  \eqref{E:subadd} implies that $nL_n\to\ell$ as $n\to\infty$,
where $\ell:=\inf_{n\ge1}nL_n$.
The question was raised in \cite{GR05}
and again in \cite[Problem~H]{BR18} as to whether
 $\ell>0$. Our second result gives a negative answer to this question
 by showing that $L_n=O(n^{-3/2})$ as $n\to\infty$. Here is the precise statement.

 \begin{theorem}\label{T:ub}
 For all $n\ge2$, we have
 \begin{equation}\label{E:ub}
L_n\le \frac{(9/2)\sqrt{\kappa}}{\bigl(\sqrt{n+\kappa}-2\sqrt{\kappa}\bigr)^2\bigl(\sqrt{n+\kappa}+\sqrt{\kappa}\bigr)},
\end{equation}
where $\kappa:=2\pi/(9\sqrt{3})$.
In particular,
\begin{equation}\label{E:limsup}
L_n\le \sqrt{\frac{3\pi\sqrt{3}}{2}}n^{-3/2}+O(n^{-2}) \quad\text{as~}n\to\infty.
\end{equation}
\end{theorem}


\section{Proof of Theorem~\ref{T:lb}}

We begin with some elementary remarks about algebroid multifunctions.
If $G,H$ are algebroid multifunctions on $D$ of degree $l,m$ respectively,
then the multifunction $(G\cup H)(z):=G(z)\cup H(z)$ is  algebroid of degree $l+m$.
The following lemma,  sometimes known as the localization principle,
gives conditions under which this process can be reversed.

\begin{lemma}[\protect{\cite[Lemma~7]{GR05}}]\label{L:localization}
Let $F$ be an algebroid multifunction on a domain $D$. Let $U,V$ be
disjoint open subsets of $\CC$ such that $F(D)\subset U\cup V$ and 
$F(D)\cap U\ne\emptyset$ and $F(D)\cap V\ne\emptyset$. 
Define $G(z):=F(z)\cap U$ and $H(z):=F(z)\cap V~(z\in D)$.
Then $G,H$ are algebroid multifunctions on $D$.
\end{lemma}

We shall also need a simple result about matching distance.
Given two multisets of $n$ points, say $\{\alpha_1,\dots,\alpha_n\}$
and $\{\beta_1,\dots,\beta_n\}$, we define the \emph{matching
distance} between them by
\[
d(\{\alpha_j\},\{\beta_j\}):=\min_{\sigma}\max_{1\le j\le n}|\alpha_j-\beta_{\sigma(j)}|,
\]
where the minimum is taken over all permutations $\sigma$
of $\{1,\dots,n\}$.

\begin{lemma}\label{L:matchingd}
Let $F$ be an algebroid multifunction of degree $n$
on a domain $D$. Let $z_0\in D$, and 
suppose that $\dist(w,F(z_0))\le \delta$ for all $w\in F(D)$.
Then
\begin{equation}\label{E:matchingd}
d(F(z),F(z_0))\le 2n\delta \quad(z\in D).
\end{equation}
\end{lemma}

\begin{proof}
Let $\delta'>\delta$. Then, for each $z\in D$, the multiset $F(z)$
is contained in the union $U$ of the open disks of radius $\delta'$
about the points of $F(z_0)$. By Lemma~\ref{L:localization},
if $V$ is a connected component of $U$, then $F(z)\cap V$
is an algebroid multifunction of degree $m$, where $1\le m\le n$.
Thus $V$ is a connected union of at most $m$ disks of radius $\delta'$, 
so has diameter 
at most $2m\delta'\le 2n\delta'$. Matching the roots of $F(z)\cap V$ with those of $F(z_0)\cap V$, and repeating for each connected component $V$
of $U$, we see that $d(F(z),F(z_0))\le 2n\delta'$.
Finally, letting $\delta'\downarrow \delta$, we obtain \eqref{E:matchingd}.
\end{proof}

We now turn to the proof of Theorem~\ref{T:lb}.
The proof proceeds via two lemmas.

\begin{lemma}\label{L:lb1}
Let $F$ be an algebroid multifunction of degree $n$ defined
on the unit disk $\DD$. Assume that its image $F(\DD)$ 
has finite diameter $M$ and inradius $\rho$.
Then
\begin{equation}\label{E:lb1}
d(F(z),F(0))\le 4n\rho^{1/(3n)}M^{1-1/(3n)} \quad(|z|\le1/2).
\end{equation}
\end{lemma}

\begin{proof}
We claim that
\begin{equation}\label{E:claim}
\dist(w, F(0))\le 2\rho^{1/(3n)}M^{1-1/(3n)} \quad(w\in F(\zeta),\,|\zeta|<1/2).
\end{equation}
If so, then  applying Lemma~\ref{L:matchingd} with $D=\{|\zeta|<1/2\}$
gives \eqref{E:lb1} for all $|z|<1/2$, and by continuity 
we deduce that it in fact holds for all $|z|\le 1/2$.

Fix $\zeta$ with $|\zeta|< 1/2$ and let $w\in F(\zeta)$. 
Let $\epsilon>0$. Since $F(\DD)$ contains no disk
of radius $\rho+\epsilon$, there exists $b\notin F(\DD)$
such that $|b-w|\le\rho+ \epsilon$. We can also suppose
that $\diam(F(\DD)\cup\{b\})\le M+\epsilon$:
indeed, we may first choose a point of the boundary of 
$F(\DD)$ at distance at most $\rho$ from $w$, and then take $b$ outside $F(\DD)$ sufficiently close to this boundary point.

Let us write
\[
F(z)=\{w\in\CC:P(z,w)=0\},
\]
where $P(z,w)=w^n+a_1(z)w^{n-1}+\cdots+a_n(z)$.
If $z\in\DD$, then $|P(z,b)|$ is the product of the distances of $b$ to the roots of $P(z,\cdot)$, 
namely the elements of $F(z)$. As $\diam (F(z)\cup\{b\})\le M+\epsilon$, it follows that
\[
|P(z,b)|\le (M+\epsilon)^n \quad(z\in\DD).
\]
Also, since $b\notin F(\DD)$, it follows that $P(z,b)\ne0$
for all $z\in\DD$. Thus $u(z):=\log((M+\epsilon)^n/|P(z,b)|)$
is a non-negative harmonic function on $\DD$.
By Harnack's inequality, $u(z)\le u(0)(1+|z|)/(1-|z|)$ for all $z\in\DD$.
In particular, $u(\zeta)\le 3u(0)$. This translates to the inequality
\begin{equation}\label{E:Harnack}
|P(0,b)|^3\le (M+\epsilon)^{2n}|P(\zeta,b)|.
\end{equation}

We shall now find estimates for each side of this inequality.
First of all, as $|b-w|\le\rho+\epsilon$,
it follows that $b$ is a distance at least 
$\dist(w,F(0))-(\rho+\epsilon)$ from $F(0)$, and thus,
as $|P(0,b)|$ is the product of the distances of $b$
to the elements of $F(0)$,
\[
|P(0,b)|\ge \bigl(\dist(w,F(0))-(\rho+\epsilon)\bigr)_+^{n},
\]
where $x_+:=\max\{x,0\}$.
Also, since $w$ is a root of $P(\zeta,\cdot)$ and
$|w-b|\le \rho+\epsilon$, a similar argument gives
\[
|P(\zeta,b)|\le (\rho+\epsilon)(M+\epsilon)^{n-1}.
\]
Inserting these inequalities into \eqref{E:Harnack}, we obtain
\[
 \bigl(\dist(w,F(0))-(\rho+\epsilon)\bigr)_+^{3n}
 \le (\rho+\epsilon)(M+\epsilon)^{3n-1}.
 \]
 Simplifying and letting $\epsilon\to0$, we deduce that
 \[
 \dist(w,F(0))\le \rho+\rho^{1/(3n)}M^{1-1/(3n)}.
 \]
As $\rho\le M$, we have $\rho\le\rho^{1/(3n)}M^{1-1/(3n)}$, and so
\[
 \dist(w,F(0))\le 2\rho^{1/(3n)}M^{1-1/(3n)}.
 \]
This justifies the claim \eqref{E:claim}.
\end{proof}

In the second lemma, we use a renormalization argument 
to remove the boundedness assumption on $F(\DD)$.

\begin{lemma}\label{L:lb2}
Let $F$ be an algebroid multifunction of degree $n$ defined
on the unit disk $\DD$. Assume that its image $F(\DD)$ 
has inradius $\rho$.
Then
\begin{equation}\label{E:lb2}
d(F(z),F(0))\le 2^{15n-3}n^{6n-1}\rho \quad(|z|\le1/7).
\end{equation}
\end{lemma}

\begin{proof}
Let $\delta>0$, and suppose that 
\begin{equation}\label{E:delta}
\sup_{|z|\le 1/7}d(F(z),F(0))>\delta.
\end{equation}
Consider all pairs of complex numbers $u,v$ 
with $|u|\le 3/4$ and $|v|\le 3/4$ and such that $d(F(u),F(v))\ge\delta$.
Among them, choose $(z_0,z_1)$ maximizing
\[
\frac{3/4-\max\{|u|,|v|\}}{|u-v|}.
\]
This maximum exists by compactness and continuity of $F$;
admissible pairs stay away from the diagonal.
Put $h:=|z_1-z_0|$ and $s:=3/4-\max\{|z_0|,|z_1|\}$.
The relation \eqref{E:delta} implies that there is an 
admissible pair $(z,0)$ with $|z|\le 1/7$, and consequently
\[
\frac{s}{h}\ge \frac{3/4-|z|}{|z|}=\frac{3}{4|z|}-1\ge 17/4>4.
\]
If $u,v$ lie in the open disk $D(z_0,2h)$
and satisfy $0<|u-v|<h/2$, then 
\[
\frac{3/4-\max\{|u|,|v|\}}{|u-v|}
\ge \frac{3/4-|z_0|-2h}{h/2}
\ge \frac{2s}{h}-4
>\frac{s}{h},
\]
and therefore, $(u,v)$ is not an admissible pair,
i.e., $d(F(u),F(v))< \delta$.
Dividing any segment from $z_0$ to a point of $D(z_0,2h)$
into four equal parts and using the triangle inequality yields
\[
d(F(z),F(z_0))<4\delta \quad(z\in D(z_0,2h)).
\]
(Note that $D(z_0,2h)\subset D(0,3/4)$ because 
$|z_0|+2h\le 3/4-s+2h < 3/4$.)

We now renormalize. Define 
\[
G(\zeta):=F(z_0+2h\zeta) \quad(\zeta\in\DD).
\]
This is an algebroid multifunction of degree $n$ satisfying
\[
d(G(\zeta),G(0))<4\delta \quad (\zeta\in\DD).
\]
In particular, $G(\DD)\subset U$, where 
$U$ is the union of the open disks of radius $4\delta$
with centres in $G(0)$. Let $V$ be a connected component
of $U$. By Lemma~\ref{L:localization}, the multifunction
$G(\zeta)\cap V$ is algebroid of degree~$m$ for some $m\le n$.
Since $V$ is a connected union of at most $m$ disks
of radius $4\delta$, its diameter is at most $8m\delta$.
In particular, $G(\DD)\cap V$ has diameter at most $8m\delta$.
Clearly its inradius is at most that of $F(\DD)$, namely~$\rho$.
Applying Lemma~\ref{L:lb1} to $G(\zeta)\cap V$, we deduce that,
whenever $|\zeta|\le 1/2$,
\begin{equation}\label{E:ineq1}
d(G(\zeta)\cap V,\,G(0)\cap V)
\le 4m\rho^{1/(3m)}(8m\delta)^{1-1/(3m)}.
\end{equation}
On the other hand, if we set $\zeta_1:=(z_1-z_0)/(2h)$, then
$|\zeta_1|=1/2$ and
\[
d(G(\zeta_1),G(0))=d(F(z_1),F(z_0))\ge \delta,
\]
so there must be at least one component $V$ of $U$
for which
\begin{equation}\label{E:ineq2}
d(G(\zeta_1)\cap V,\,G(0)\cap V)\ge \delta.
\end{equation}
Combining the inequalities \eqref{E:ineq1} and \eqref{E:ineq2}
gives
\[
\delta\le 4m\rho^{1/(3m)}(8m\delta)^{1-1/(3m)},
\]
which, after simplification, becomes
\[
\delta\le 2^{15m-3}m^{6m-1}\rho.
\]
As $m\le n$, it follows that
\[
\delta\le 2^{15n-3}n^{6n-1}\rho.
\]
Finally, as this holds for each $\delta>0$ satisfying \eqref{E:delta},
the inequality \eqref{E:lb2} follows.
\end{proof}

\begin{proof}[Proof of Theorem~\ref{T:lb}]
Let $F$ be an algebroid multifunction of degree $n$
on $\DD$. By Lemma~\ref{L:lb2}, it satisfies
\begin{equation}\label{E:complicatedineq}
d(F(z),F(0))\le 2^{15n-3}n^{6n-1}\rho(F(\DD))
\quad(|z|\le 1/7).
\end{equation}
Writing $F$ as in \eqref{E:algdef}, the function $-a_1(z)$
is the sum of the roots, i.e., the sum of the elements
of $F(z)$, counted according to multiplicity. It follows 
from \eqref{E:complicatedineq} that
\begin{equation}\label{E:ultimate}
|a_1(z)-a_1(0)|\le n\cdot2^{15n-3}n^{6n-1}\rho(F(\DD))
\quad(|z|\le 1/7).
\end{equation}
On the other hand, if $F$ is normalized so that $a_1'(0)=1$,
then, by the standard Cauchy estimate,
\begin{equation}\label{E:Cauchy}
\sup_{|z|=1/7}\frac{|a_1(z)-a_1(0)|}{1/7}\ge |a_1'(0)|=1.
\end{equation}
Combining \eqref{E:ultimate} and \eqref{E:Cauchy},
we  deduce that
\[
\rho(F(\DD))\ge \frac{8/7}{2^{15n}n^{6n}}.
\]
As this estimate holds for each $F$ of degree $n$
for which $a_1'(0)=1$, we conclude that
\[
L_n\ge \frac{8/7}{2^{15n}n^{6n}}>\frac{1}{2^{15n}n^{6n}}.\qedhere
\]
\end{proof}


\section{Proof of Theorem~\ref{T:ub}}

The proof of Theorem~\ref{T:ub} also proceeds via 
two lemmas. We write $|E|$ for the cardinality
of a finite set $E$.

\begin{lemma}\label{L:ub1}
Let $E$ be a finite, non-empty subset of $\DD$.
Then, for each $n>|E|$, we have
\begin{equation}\label{E:ub1}
L_n\le \frac{\rho(\DD\setminus E)}{n-|E|}.
\end{equation}
\end{lemma}

\begin{proof}
Let $n>|E|$.
Fix $0<t<1$ and  set 
\[
m:=n-|E|
\quad\text{and}\quad
\epsilon:=\frac{1}{3}\dist(E,\partial\DD)^{n-m}(1-t)^m.
\]
Define
\[
F(z):=\Bigl\{w\in\CC: \prod_{c\in E}(w-c)\cdot(w+tz)^m+\epsilon=0\Bigr\}
\quad(z\in \DD).
\]
Then $F$ is an algebroid multifunction of degree~$n$ on $\DD$.

We first show that $F(\DD)\subset\DD$.
Fix $z\in\DD$. If $w\in\partial\DD$, then
\[
\Bigl| \prod_{c\in E}(w-c)\cdot(w+tz)^m+\epsilon\Bigr|
\ge \dist(E,\partial\DD)^{n-m}(1-t)^m-\epsilon>\epsilon.
\]
Therefore, by Rouch\'e's theorem, the number of roots in $\DD$
of $w\mapsto\prod_{c\in E}(w-c)\cdot(w+tz)^m+\epsilon$ 
is the same as the number of roots in $\DD$ of  
$w\mapsto\prod_{c\in E}(w-c)\cdot(w+tz)^m$,
namely~$n$. Hence $F(z)\subset\DD$, as desired.

We next show that $F(\DD)\cap E=\emptyset$.
This follows immediately from the fact that, if $z\in\DD$ and
$w\in E$, then
\[
\prod_{c\in E}(w-c)\cdot(w+tz)^m+\epsilon=\epsilon\ne0.
\]
Combining this observation with the inclusion $F(\DD)\subset\DD$,
we obtain $F(\DD)\subset\DD\setminus E$, and hence
\[
\rho(F(\DD))\le \rho(\DD\setminus E).
\]

Now we consider the $a_1$ function associated to $F$.
The coefficient
of $w^{n-1}$ in $\prod_{c\in E}(w-c)\cdot(w+tz)^m+\epsilon$
is  $a_1(z)=mtz-\sum_{c\in E}c$. Thus $a_1'(0)=mt$.
Setting $\tilde{F}:=(1/mt)F$, we obtain a new algebroid
multifunction of degree $n$ whose $a_1$ function satisfies
$a_1'(0)=1$ and whose image has inradius satisfying
\[
\rho(\tilde{F}(\DD))=\frac{\rho(F(\DD))}{mt}\le \frac{\rho(\DD\setminus E)}{mt}.
\]
Hence
\[
L_n\le  \frac{\rho(\DD\setminus E)}{mt}.
\]
Letting $t\to1$, we get
\[
L_n\le \frac{\rho(\DD\setminus E)}{m}=\frac{\rho(\DD\setminus E)}{n-|E|}
\]
which is just \eqref{E:ub1}.
\end{proof}

The next lemma furnishes an efficient choice of $E$.
This result is standard, but a short proof is included
for completeness.

\begin{lemma}\label{L:ub2}
Given $r>0$, there exists a non-empty, finite subset $E_r$ of $\DD$
such that 
\begin{equation}\label{E:ub2}
\rho(\DD\setminus E_r)\le r
\quad\text{and}\quad
|E_r|\le \frac{2\pi}{3\sqrt{3}}\frac{(1+r)^2}{r^2}.
\end{equation}
\end{lemma}

\begin{proof}
Let $\Lambda_r$ be the triangular lattice given by
\[
\Lambda_r:=\{\sqrt{3}r(m+ne^{i\pi/3}):m,n\in\ZZ\}.
\]
Its nearest-neighbour separation is $\sqrt{3}r$, 
and each Voronoi cell is a regular hexagon with
circumradius $r$.
Set $E_r:=\Lambda_r\cap \DD$.

We claim that 
\begin{equation}\label{E:rhoest}
\rho(\DD\setminus E_r)\le r.
\end{equation}
Indeed, suppose the contrary.
Then there exists an open subdisk $D$ of $\DD\setminus E_r$
of radius $s>r$. Let $c$ be the centre of $D$.
There exists a point $b\in\Lambda_r$ with $|b-c|\le r$.
This implies that $b\in D$.
Hence $b\in D\cap\Lambda_r=D\cap\DD\cap\Lambda_r=D\cap E_r$,
contrary to the choice of $D$. This justifies the claim
\eqref{E:rhoest}.

It remains to show that
\begin{equation}\label{E:areaest}
|E_r|\le \frac{2\pi}{3\sqrt{3}}\frac{(1+r)^2}{r^2}.
\end{equation}
For this, we remark that each Voronoi cell 
with centre in $E_r$ lies in the closed disk $\overline{D}(0,1+r)$.
These hexagonal cells each have area  $3\sqrt{3}r^2/2$ and are disjoint up to boundaries, so a comparison of areas gives
\[
|E_r|\frac{3\sqrt{3}r^2}{2}\le \pi(1+r)^2.
\]
This establishes \eqref{E:areaest}.
\end{proof}

\begin{proof}[Proof of Theorem~\ref{T:ub}]
Let $n\ge2$ and $\kappa:=2\pi/(9\sqrt{3})$.
By Lemma~\ref{L:ub2}, for each $r>0$, there exists a subset $E_r$ of $\DD$ 
such that $\rho(\DD\setminus E_r)\le r$ and $|E_r|\le3\kappa(1+r)^2/r^2$.
Applying Lemma~\ref{L:ub1} with $E=E_r$, we deduce that $L_n\le \rho(\DD\setminus E_r)/(n-|E_r|)$
for all $r$ such that $|E_r|<n$. It follows that
\[
L_n\le \frac{r}{n-3\kappa(1+r)^2/r^2}
\]
for all $r>0$ such that the denominator is positive.
In particular, if we take 
\[
r=\frac{3\sqrt{\kappa}}{\sqrt{n+\kappa}-2\sqrt{\kappa}},
\]
then 
\[
n-3\kappa\frac{(1+r)^2}{r^2}=\frac{2}{3}\bigl(\sqrt{n+\kappa}-2\sqrt{\kappa}\bigr)\bigl(\sqrt{n+\kappa}+\sqrt{\kappa})>0,
\]
and hence
\[
L_n\le \frac{(9/2)\sqrt{\kappa}}{\bigl(\sqrt{n+\kappa}-2\sqrt{\kappa}\bigr)^2\bigl(\sqrt{n+\kappa}+\sqrt{\kappa}\bigr)}.
\]
This proves \eqref{E:ub}, and the asymptotic bound \eqref{E:limsup} is an easy consequence.
\end{proof}


\section{Concluding remarks}

Evidently, there is a huge gap between the lower and upper bounds \eqref{E:lb} and  \eqref{E:ub}.
We suspect that the true value of $L_n$ is closer to the upper bound. 
In fact, we would even go so far
as to guess that
\[
\liminf_{n\to\infty}n^{3/2}L_n>0.
\]
However, for the moment, a proof of  this seems out of reach.


\section*{AI tools declaration}
In preparing this work, the author used ChatGPT Pro Astra to develop and explore strategies for the proofs of 
Theorems~\ref{T:lb} and \ref{T:ub}.
The final mathematical  text was written entirely by the author, who thoroughly corrected, 
verified, and streamlined the proposed strategy. 
The author accepts full responsibility for the accuracy and integrity of the final results presented in this paper.

\bibliographystyle{plain}
\bibliography{bibfile}

\end{document}